\def\version{0.28}

\def\journal{JFA}
\def\titlep{Inductive limit violates quasi-cocommutativity}
\documentclass[11pt]{article}
\usepackage{graphicx,ifthen}
\usepackage{amssymb}
\usepackage{amsmath}

\font\germ=eufm10 at12pt

\def\goth#1{\hbox{\germ#1}}

\newcommand{\qed}{\hbox{\rule[-2pt]{3pt}{6pt}}}
\newcommand{\qedh}{\hfill\qed \\}

\newcommand{\vv}{\vspace{.3in}}

\newcommand{\vep}{\varepsilon}

\newtheorem{Thm}{Theorem}[section]

\newtheorem{rem}[Thm]{Remark}

\newtheorem{defi}[Thm]{Definition}
\newtheorem{lem}[Thm]{Lemma}

\newtheorem{cor}[Thm]{Corollary}

\newcommand{\kn}{\Large\bf
$K\hspace{-.4cm} N$
\Large\bf\vv }

\def\cal#1{\mathcal #1}

\def\pr{{\it Proof.}\quad}

\def\co#1{{\cal O}_{#1}}
\def\disp#1{{\displaystyle #1}}
\def\brl{branching law}
\def\bfsnl{{\rm BFS}_{N}(\Lambda)}

\addtocounter{footnote}{1}
\def\sftt#1{
\setcounter{equation}{0}
\addtocounter{footnote}{1}
\section{#1}
}

\def\ssft#1{\subsection{#1}}

\def\ngt{{\bf N}_{\geq 2}}

\begin{document}
%
%
\def\autherp{Katsunori Kawamura}
\def\emailp{e-mail: kawamura@kurims.kyoto-u.ac.jp.}
\def\addressp{{\small {\it College of Science and Engineering,
 Ritsumeikan University,}}\\
{\small {\it 1-1-1 Noji Higashi, Kusatsu, Shiga 525-8577, Japan}}
}

\def\infw{\Lambda^{\frac{\infty}{2}}V}
\def\zhalfs{{\bf Z}+\frac{1}{2}}
\def\ems{\emptyset}
\def\pmvac{|{\rm vac}\!\!>\!\! _{\pm}}
\def\vac{|{\rm vac}\rangle _{+}}
\def\dvac{|{\rm vac}\rangle _{-}}
\def\ovac{|0\rangle}
\def\tovac{|\tilde{0}\rangle}
\def\expt#1{\langle #1\rangle}
\def\zph{{\bf Z}_{+/2}}
\def\zmh{{\bf Z}_{-/2}}
\def\brl{branching law}
\def\bfsnl{{\rm BFS}_{N}(\Lambda)}
\def\scm#1{S({\bf C}^{N})^{\otimes #1}}
\def\mqb{\{(M_{i},q_{i},B_{i})\}_{i=1}^{N}}
\def\zhalf{\mbox{${\bf Z}+\frac{1}{2}$}}
\def\zmha{\mbox{${\bf Z}_{\leq 0}-\frac{1}{2}$}}
\newcommand{\mline}{\noindent
\thicklines
\setlength{\unitlength}{.1mm}
\begin{picture}(1000,5)
\put(0,0){\line(1,0){1250}}
\end{picture}
\par
 }
\def\ptimes{\otimes_{\varphi}}
\def\delp{\Delta_{\varphi}}
\def\delps{\Delta_{\varphi^{*}}}
\def\gamp{\Gamma_{\varphi}}
\def\gamps{\Gamma_{\varphi^{*}}}
\def\sem{{\sf M}}
\def\hdelp{\hat{\Delta}_{\varphi}}
\def\tilco#1{\tilde{\co{#1}}}
\def\ndm#1{{\bf M}_{#1}(\{0,1\})}
\def\cdm#1{{\cal M}_{#1}(\{0,1\})}
\def\tndm#1{\tilde{{\bf M}}_{#1}(\{0,1\})}
\def\sck{{\sf CK}_{*}}
\def\hdel{\hat{\Delta}}
\def\ba{\mbox{\boldmath$a$}}
\def\bb{\mbox{\boldmath$b$}}
\def\bc{\mbox{\boldmath$c$}}
\def\bd{\mbox{\boldmath$d$}}
\def\be{\mbox{\boldmath$e$}}
\def\bk{\mbox{\boldmath$k$}}
\def\bm{\mbox{\boldmath$m$}}
\def\bp{\mbox{\boldmath$p$}}
\def\bq{\mbox{\boldmath$q$}}
\def\bu{\mbox{\boldmath$u$}}
\def\bv{\mbox{\boldmath$v$}}
\def\bw{\mbox{\boldmath$w$}}
\def\bx{\mbox{\boldmath$x$}}
\def\by{\mbox{\boldmath$y$}}
\def\bz{\mbox{\boldmath$z$}}
\def\bomega{\mbox{\boldmath$\omega$}}
\def\N{{\bf N}}
\def\lxm{L_{2}(X,\mu)}
\def\qtimes{\otimes_{\tilde{\varphi}}}
\def\ul#1{\underline{#1}}
\def\ngt{{\bf N}_{\geq 2}}
\def\ngti{{\bf N}_{\geq 2}^{\infty}}
\def\mninf{M_{n}^{\otimes \infty}}
\def\titlepage{

\noindent
{\bf 
\noindent
\thicklines
\setlength{\unitlength}{.1mm}
\begin{picture}(1000,0)(0,-300)
\put(0,0){\kn \knn\, for \journal\, Ver.\,\version}
\put(0,-50){\today,\quad {\rm file:} {\rm {\small \textsf{tit01.txt,\, J1.tex}}}}
\end{picture}
}
\vspace{-2.5cm}
\quad\\
{\small 
\footnote{
${\displaystyle
\mbox{directory: \textsf{\fileplace}, 
file: \textsf{\incfile},\, from \startdate}}$}}
\quad\\
\framebox{
\begin{tabular}{ll}
\textsf{Title:} &
\begin{minipage}[t]{4in}
\titlep
\end{minipage}
\\
\textsf{Author:} &\autherp
\end{tabular}
}
{\footnotesize	
\tableofcontents }
}

%
%
%
\setcounter{section}{0}
\setcounter{footnote}{0}
\setcounter{page}{1}
\pagestyle{plain}

%
%
\title{\titlep}
\author{\autherp\thanks{\emailp}
\\
\addressp}
\date{}
\maketitle

%
%
\begin{abstract}
We show that the inductive limit of a certain inductive system 
of quasi-cocommutative C$^{*}$-bialgebras is not quasi-cocommutative.
This implies that
the category of quasi-cocommutative C$^{*}$-bialgebras
is not closed with respect to the inductive limit. 
\end{abstract}


%
%
\sftt{Introduction}
\label{section:first}
A C$^{*}$-bialgebra is a C$^{*}$-algebra
with several extra structures on it,
which  was introduced in
the operator algebra approach to quantum groups 
as a locally compact quantum semigroup
\cite{KV} (see also \cite{ES,MNW}).
A quasi-cocommutative C$^{*}$-bialgebra is defined as a C$^{*}$-bialgebra
with a universal $R$-matrix which is 
modified to focus on C$^{*}$-bialgebra
\cite{Kassel,TS21}.
In this paper,
we prove the following statement. 
%
%
\begin{Thm}
\label{Thm:mainone}
The category of quasi-cocommutative C$^{*}$-bialgebras
is not closed with respect to the inductive limit.
\end{Thm}
In this section, we roughly explain our motivation 
and the significance of Theorem \ref{Thm:mainone}.
Explicit mathematical definitions will 
be shown after $\S$ \ref{section:second}.
%
%
\ssft{Motivation}
\label{subsection:firstone}
We have constructed various non-commutative and non-cocommutative 
C$^{*}$-bialgebras by using sets of C$^{*}$-algebras
and $*$-homomorphisms among them \cite{TS02,TS05,TS20}.
These studies are motivated by 
a discovery of a certain non-symmetric tensor product 
of representations of Cuntz algebras \cite{TS01},
without reference to the study of quantum groups.

On the other hand, the inductive limit is an important 
tool to construct new C$^{*}$-algebras, which 
often changes properties of  C$^{*}$-algebras \cite{Glimm,Powers}.
For a given subcategory of the category of C$^{*}$-algebras, 
the inductive limit on it is an interesting subject of research.

Our interests are to define the inductive limit of C$^{*}$-bialgebras
and to study its property.
Especially, we consider the inductive limit of quasi-cocommutative
C$^{*}$-bialgebras in this paper.

In order to prove Theorem \ref{Thm:mainone},
we construct an inductive system of
quasi-cocommutative C$^{*}$-bialgebras
such that its inductive limit is not quasi-cocommutative.

%
%
\ssft{Comparison with quantum enveloping algebras}
\label{subsection:firsttwo}
We explain the significance of Theorem \ref{Thm:mainone}
in comparison to cases of quantum enveloping algebras in this subsection.

We start with a brief history of quantum groups.
At the beginning,
quantum groups were introduced as one-parameter deformations of universal
enveloping algebras $U_{q}({\goth g})$ 
(= {\it quantum enveloping algebra} \cite{Kassel})
of semisimple complex Lie algebras ${\goth g}$ \cite{Drinfeld,Jimbo}.
A motivation of the study is to construct solutions 
of Yang-Baxter equations \cite{Baxter,Yang},
which are called $R$-matrices.
The fundamental structure of a quantum group is a Hopf algebra,
but it is not sufficient for the original purpose.
In order to define its universal $R$-matrix as an infinite series
in the tensor square of the completion of $U_{q}({\goth g})$,
the $h$-adic topology is introduced (\cite{Kassel}, Chap. XVI)
where $h$ is related to $q$ as $q=e^{h}$.
In this case, the topology is an inverse limit topology,
and 
it is necessary to define the universal $R$-matrix in 
the theory of quantum groups in general.
The quasi-cocommutativity is acquired by taking the limit in this case.

On the other hand,
Theorem \ref{Thm:mainone} means that inductive limit violates 
the quasi-cocommutativity of C$^{*}$-bialgebras in general.
This is a new phenomenon which is quite a contrast 
to the case of quantum enveloping algebras.
Both cases give an indication of relations between
topology and quasi-cocommutativity.

In $\S$ \ref{section:second},
we will recall basic definitions of 
quasi-cocommutative C$^{*}$-bialgebras and consider categories.
In $\S$ \ref{section:third},
we will prepare a general method to construct quasi-cocommutative C$^{*}$-bialgebras.
In $\S$ \ref{section:fourth},
we will construct an example of inductive system of quasi-cocommutative
C$^{*}$-bialgebras and prove Theorem \ref{Thm:mainone}.

%
%
\sftt{Definitions}
\label{section:second}
In this section, we recall basic definitions.
%
%
\ssft{Quasi-cocommutative C$^{*}$-bialgebra and categories}
\label{subsection:secondone}
In this subsection,
we recall definitions of quasi-cocommutative C$^{*}$-bialgebra 
and related notions \cite{TS20,TS21}.
In addition,
we consider categories of C$^{*}$-bialgebras in Remark \ref{rem:cat}.

For a C$^{*}$-algebra $A$,
let ${\cal M}(A)$  denote the  multiplier algebra of $A$. 
For two C$^{*}$-algebras $A$ and $B$,
let ${\rm Hom}(A,B)$ and $A\otimes B$ 
denote the set of all $*$-homomorphisms from $A$ to $B$ and 
the minimal C$^{*}$-tensor product of $A$ and $B$, respectively.
We state that $f\in {\rm Hom}(A,{\cal M}(B))$ 
is {\it nondegenerate} if $f (A)B$ is dense in $B$.
In this case,
$f$ is called a {\it morphism} from $A$ to $B$ \cite{W3}.
If $f$ is a nondegenerate $*$-homomorphism from $A$ to $B$,
then we can regard $f$ as a morphism from $A$ to $B$
by using the canonical embedding of $B$ into ${\cal M}(B)$.
Each morphism $f$ from $A$ to $B$ can be extended uniquely 
to a homomorphism $\tilde{f}$ from ${\cal M}(A)$ to ${\cal M}(B)$
such that $\tilde{f}(m) f(b)a = f(mb)a$ 
for $m\in {\cal M}(B), b\in B$, and $a \in A$.

A {\it C$^{*}$-bialgebra} is a pair $(A,\Delta)$
of a C$^{*}$-algebra $A$ and a morphism $\Delta$
from $A$ to $A\otimes A$ which satisfies
$(\Delta\otimes id)\circ \Delta=(id\otimes\Delta)\circ \Delta$.
We call $\Delta$ the {\it comultiplication} of $A$.
A C$^{*}$-bialgebra $(A,\Delta)$ is {\it counital}
if there exists $\vep\in {\rm Hom}(A,{\bf C})$ such that
$(\vep\otimes id)\circ \Delta= id = (id\otimes \vep)\circ
\Delta$.
We call $\vep$ the {\it counit} of $A$ and write $(A,\Delta,\vep)$ as 
the counital C$^{*}$-bialgebra $(A,\Delta)$ with the counit $\vep$.
Remark that we do not assume $\Delta(A)\subset A\otimes A$.
Furthermore,
$A$ has no unit for a C$^{*}$-bialgebra $(A,\Delta)$ in general.

Define the {\it extended flip} $\tilde{\tau}_{A,A}$ from ${\cal M}(A\otimes A)$ 
to ${\cal M}(A\otimes A)$ as
$\tilde{\tau}_{A,A}(X)(x\otimes y)\equiv \tau_{A,A}(X(y\otimes x))$
for $X\in {\cal M}(A\otimes A),\,x, y\in A$
where $\tau_{A,A}$ denotes the flip of $A\otimes A$.
The map $\Delta^{op}$ from $A$ to ${\cal M}(A\otimes A)$
defined as $\Delta^{op}\equiv \tilde{\tau}_{A,A}\circ \Delta$
is called the {\it opposite comultiplication} of $\Delta$.
A C$^{*}$-bialgebra $(A,\Delta)$ is {\it cocommutative}
if $\Delta=\Delta^{op}$. 
An element $R$ in ${\cal M}(A\otimes A)$
is called a {\it (unitary) universal $R$-matrix} of $(A,\Delta)$
if $R$ is a unitary and 
%
%
\begin{equation}
\label{eqn:univ}
R\Delta(x)R^{*}=\Delta^{op}(x)\quad
(x\in A).
\end{equation}
In this case,
we state that $(A,\Delta)$ is {\it quasi-cocommutative
(or almost cocommutative \cite{CP}).}
We write $(A,\Delta,R)$ as  
a quasi-cocommutative C$^{*}$-bialgebra $(A,\Delta)$ with
a universal $R$-matrix $R$. 
If $A$ is unital,
then ${\cal M}(A\otimes A)=A\otimes A$ and
$\tilde{\tau}_{A,A}=\tau_{A,A}$.
In addition, if $(A,\Delta)$ is quasi-cocommutative with
a universal $R$-matrix $R$, then $R\in A\otimes A$.
We state that 
a quasi-cocommutative C$^{*}$-bialgebra $(A,\Delta, R)$ 
is {\it quasi-triangular  (or braided \cite{Kassel})} if the following holds:
%
%
\begin{equation}
\label{eqn:delone}
(\Delta\otimes id)(R)=R_{13}R_{23},
\quad
(id\otimes \Delta)(R)=R_{13}R_{12}
\end{equation}
where we use the leg numbering notation \cite{BS};
$(A,\Delta, R)$ is {\it triangular} if 
$(A,\Delta, R)$ is quasi-triangular and 
the following holds:
%
%
\begin{equation}
\label{eqn:tri}
R\,\tilde{\tau}_{A,A}(R)=I
\end{equation}
where  $I$ denotes the unit of ${\cal M}(A\otimes A)$.
The cocommutativity is the dual notion of the commutativity.
Since a cocommutative C$^{*}$-bialgebra
is always quasi-cocommutative (furthermore, it is triangular), 
the quasi-cocommutativity is a generalization of 
the cocommutativity.

We consider categories of C$^{*}$-bialgebras as follows.
%
%
\begin{rem}
\label{rem:cat}
{\rm
Let $(A_{i},\Delta_{i})$ be a C$^{*}$-bialgebra for $i=1,2$.
\begin{enumerate}
\item
A map $f$ is a {\it C$^{*}$-bialgebra morphism} from 
$(A_{1},\Delta_{1})$ to $(A_{2},\Delta_{2})$ 
if $f$ is a nondegenerate $*$-homomorphism from $A_{1}$ to $M(A_{2})$
such that $(f\otimes f)\circ \Delta_{1}=\Delta_{2}\circ f$.
The category of C$^{*}$-bialgebras is defined as the
category with C$^{*}$-bialgebra morphisms as morphisms among objects.
In addition, if $R_1$ is a universal $R$-matrix of $(A_{1},\Delta_1)$,
then $R^{'}\equiv (f\otimes f)(R_1)$
is also a universal $R$-matrix of the C$^{*}$-bialgebra 
$(f(A_{1}),\Delta_{2}|_{f(A_1)})$.
\item
Define $C\equiv A_1\otimes A_2$ and 
$\Delta\equiv (id\otimes \tau\otimes id)\circ (\Delta_{1}\otimes \Delta_{2})$
where $\tau$ denotes the flip from $A_1\otimes A_2$
to $A_2\otimes A_1$ and 
$id\otimes \tau\otimes id$ is extended on ${\cal M}(A_{1}\otimes A_{1})\otimes 
{\cal M}(A_{2}\otimes A_{2})$.
Then $C$ is a C$^{*}$-bialgebra with the comultiplication $\Delta$.
We see that the tensor product of two C$^*$-bialgebra morphisms
is also a C$^*$-bialgebra morphism.
In addition, if $R_{i}$ is a universal $R$-matrix of $A_{i}$ for $i=1,2$,
then $R\equiv (id\otimes \tau\otimes id)(R_{1}\otimes R_{2})$
is also a universal $R$-matrix of $(C,\Delta)$.
From this,
we see that the tensor product of quasi-cocommutative
C$^{*}$-bialgebras is also quasi-cocommutative.
Furthermore,
we can verify that the tensor product of quasi-triangular
({\it resp}.  triangular) 
C$^{*}$-bialgebras is also quasi-triangular
({\it resp}. triangular). 
In this way,
 the category of (quasi-cocommutative,
quasi-triangular, triangular) C$^{*}$-bialgebras 
is closed with respect to the tensor product.
\end{enumerate}
}
\end{rem}
%
%
\ssft{Direct product, direct sum and inductive limit of C$^{*}$-algebras}
\label{subsection:secondtwo}
We recall direct product, direct sum and inductive limit of C$^{*}$-algebras
\cite{Blackadar2006}.
For an infinite set $\{A_{i}:i\in \Omega\}$ of C$^{*}$-algebras,
we define two C$^{*}$-algebras $\prod_{i\in\Omega} A_{i}$
and $\bigoplus_{i\in\Omega} A_{i}$ as follows:
%
%
\begin{eqnarray}
\label{eqn:prod}
\disp{\prod_{i\in\Omega} A_{i}}\equiv &
\disp{\{(a_{i}):\|(a_{i})\|\equiv \sup_{i}\|a_{i}\|<\infty\},}\\
\nonumber 
\\
\label{eqn:plus}
\disp{\bigoplus_{i\in\Omega} A_{i}}\equiv &
\disp{\{(a_{i}):\|(a_{i})\|\to 0\mbox{ as }i\to\infty\}}
\end{eqnarray}
in the sense that for every $\vep>0$ there are only finitely many
$i$ for which $\|a_{i}\|>\vep$.
We call $\prod_{i\in\Omega} A_{i}$ and $\bigoplus_{i\in\Omega} A_{i}$ 
the {\it direct product} and the {\it direct sum} of $A_{i}$'s,
respectively. 
The algebra $\bigoplus_{i\in\Omega} A_{i}$ 
is a closed two-sided ideal of $\prod_{i\in\Omega} A_{i}$.
The algebraic direct sum $\oplus_{alg}\{A_{i}:i\in\Omega\}$ 
is a dense $*$-subalgebra of $\oplus\{A_{i}:i\in\Omega\}$.
Since ${\cal M}(\oplus_{i\in \Omega}A_{i})
\cong \prod_{i\in \Omega}{\cal M}(A_{i})$ (\cite{Blackadar2006}, II.8.1.3),
if $A_{i}$ is unital for each $i$, then
%
%
\begin{equation}
\label{eqn:multiplus}
{\cal M}\Bigl(\bigoplus_{i\in \Omega}A_{i}\Bigr)\cong \prod_{i\in \Omega}A_{i}.
\end{equation}

An {\it inductive system} of C$^{*}$-algebras is a collection
$\{(A_{i},f_{ij}):i,j\in \Omega,\, i\leq j\}$,
where $\Omega$ is a directed set, the $A_{i}$ are
C$^{*}$-algebras, and $f_{ij}$ is a $*$-homomorphism
from $A_i$ to $A_j$ with $f_{ik}=f_{jk}\circ f_{ij}$
for $i\leq j\leq k$. 
With respect to the seminorm
$\|a\|\equiv \lim_{j>i}\|f_{ij}(a)\|$ for $a\in A_{i}$,
the completion of the algebraic direct limit 
with elements of seminorm $0$ divided out
is a C$^{*}$-algebra called the {\it inductive limit}
of the system, written $\varinjlim(A_{i},f_{ij})$.
Clearly, if $A_{i}$ is commutative for each $i$,
then $\varinjlim(A_{i},f_{ij})$ is also commutative.

We introduce the inductive limit of C$^{*}$-bialgebras as follows.
%
%
\begin{defi}
\label{defi:ind}
A data $\{(A_{i},\Delta_{i}, f_{ij}):i,j\in\Omega\}$ 
is an inductive system of C$^{*}$-bialgebras 
if 
$\{(A_{i},f_{ij}):i,j\in\Omega\}$ is an inductive system 
of C$^{*}$-algebras such that 
$A_{i}$ is a C$^{*}$-bialgebra 
and $f_{ij}$ is a C$^{*}$-bialgebra morphism from $A_{i}$ to $A_{j}$.
\end{defi}

For an inductive system 
$\{(A_{i},\Delta_{i},f_{ij}):i,j\in\Omega\}$ of C$^{*}$-bialgebras,
let $A$ denote the inductive limit of the inductive system 
$\{(A_{i},f_{ij}):i,j\in\Omega\}$ of C$^{*}$-algebras.
Let $\mu_{i}$ denote the canonical map from $A_{i}$ to $A$.
Define the map $\Delta^{(0)}$ on $\bigcup_{i}\mu_{i}(A_{i})$ as 
%
%
\begin{equation}
\label{eqn:zero}
\Delta^{(0)}(\mu_{i}(x))
\equiv \{(\mu_{i}\otimes \mu_{i})\circ \Delta_{i}\}(x)\quad (x\in A_{i}).
\end{equation}
Let $\Delta$ denote  the unique extension of 
$\Delta^{(0)}$ on $A$.
Then $(A,\Delta)$ is a C$^{*}$-bialgebra.
We call $(A,\Delta)$ the {\it inductive limit} of
$\{(A_{i},\Delta_{i},f_{ij}):i,j\in\Omega\}$.
If $(A_{i},\Delta_{i})$ is cocommutative for each $i$,
then the inductive limit $(A,\Delta)$ is also cocommutative.
In this way,
the inductive limit preserves
both the commutativity and the cocommutativity. 

We prepare a lemma  for the proof of Theorem \ref{Thm:mainone} 
in $\S$ \ref{subsection:fourthtwo}.
%
%
\begin{lem}
\label{lem:maintwob}(\cite{TS20}, Lemma 2.1)
Let $(A,\Delta)$ be a C$^*$-bialgebra.
If $(A,\Delta)$ is quasi-cocommutative,
then 
for any two nondegenerate representations
$\pi_{1}$ and $\pi_{2}$ of the C$^{*}$-algebra $A$,
$(\pi_{1}\otimes \pi_{2})\circ \Delta$ 
and
$(\pi_{2}\otimes \pi_{1})\circ \Delta$ 
are unitarily equivalent
where we write the extension of
$\pi_{i}\otimes \pi_{j}$ on ${\cal M}(A\otimes A)$
as the same notation $\pi_{i}\otimes \pi_{j}$ for $i,j=1,2$.
\end{lem}

%
%
\sftt{C$^{*}$-weakly coassociative system}
\label{section:third}
In this section, 
we recall a general method to construct C$^{*}$-bialgebras
and develop it.
%
%
\ssft{Definition}
\label{subsection:thirdone}
According to \cite{TS02,TS21},
we recall C$^{*}$-weakly coassociative system in this subsection.
We call $\sem$ a {\it monoid} if $\sem$ is a semigroup with unit.
%
%
\begin{defi}
\label{defi:axiom}
Let $\sem$ be a monoid with the unit $e$.
A data $\{(A_{a},\varphi_{a,b}):a,b\in \sem\}$
is a C$^{*}$-weakly coassociative system (= C$^{*}$-WCS) over $\sem$ if 
$A_{a}$ is a unital C$^{*}$-algebra for $a\in \sem$
and $\varphi_{a,b}$ is a unital $*$-homomorphism
from $A_{ab}$ to $A_{a}\otimes A_{b}$
for $a,b\in \sem$ such that
\begin{enumerate}
\item
for all $a,b,c\in \sem$, the following holds:
%
%
\begin{equation}
\label{eqn:wcs}
(id_{a}\otimes \varphi_{b,c})\circ \varphi_{a,bc}
=(\varphi_{a,b}\otimes id_{c})\circ \varphi_{ab,c}
\end{equation}
where $id_{x}$ denotes the identity map on $A_{x}$ for $x=a,c$,
\item
there exists a counit $\vep_{e}$ of $A_{e}$ 
such that $(A_{e},\varphi_{e,e},\vep_{e})$ 
is a counital C$^{*}$-bialgebra,
\item
$\varphi_{e,a}(x)=I_{e}\otimes x$ and
$\varphi_{a,e}(x)=x\otimes I_{e}$ for $x\in A_{a}$ and $a\in \sem$.
\end{enumerate}
\end{defi}

\noindent 
Then the following holds.
%
%
\begin{Thm}
\label{Thm:mainthree}
(\cite{TS02}, Theorem 3.1)	
Let $\{(A_{a},\varphi_{a,b}):a,b\in \sem\}$ be a C$^{*}$-WCS 
over a monoid $\sem$.
Assume that $\sem$ satisfies 
%
%
\begin{equation}
\label{eqn:finiteness}
\#{\cal N}_{a}<\infty \mbox{ for each }a\in \sem
\end{equation}
where ${\cal N}_{a}\equiv\{(b,c)\in \sem\times \sem:\,bc=a\}$.
Define C$^{*}$-algebras 
%
%
\begin{equation}
\label{eqn:astar}
A_{*}\equiv  \oplus \{A_{a}:a\in \sem\},\quad
C_{a}\equiv 
\oplus \{A_{b}\otimes A_{c}:(b,c)\in {\cal N}_{a}\}
\quad (a\in\sem),
\end{equation}
and define $*$-homomorphisms $\Delta^{(a)}_{\varphi}\in{\rm Hom}(A_{a},C_{a})$ and
$\Delta_{\varphi}
\in {\rm Hom}(A_{*}, A_{*}\otimes A_{*})$ by
%
%
\begin{equation}
\label{eqn:cua}
\Delta^{(a)}_{\varphi}(x)\equiv \sum_{(b,c)\in {\cal N}_{a}}
\varphi_{b,c}(x)\quad(x\in A_{a}),\quad 
\Delta_{\varphi}\equiv \oplus\{\Delta_{\varphi}^{(a)}:a\in \sem\}.
\end{equation}
Then $(A_{*},\delp)$ is a C$^{*}$-bialgebra.
\end{Thm}

\noindent
We call $(A_{*},\Delta_{\varphi})$ in Theorem \ref{Thm:mainthree} 
the {\it C$^{*}$-bialgebra associated with} 
$\{(A_{a},\varphi_{a,b}):a,b\in \sem\}$.
In this paper, we always assume the condition (\ref{eqn:finiteness}).

%
%
\begin{defi}
\label{defi:qqq}
Let 
$\{(A_{a},\varphi_{a,b}):a,b\in\sem\}$
be a  C$^{*}$-WCS.
\begin{enumerate}
\item
For $a,b\in\sem$,
define $\varphi_{a,b}^{op}\in {\rm Hom}(A_{ab},A_{b}\otimes A_{a})$ by
%
%
\begin{equation}
\label{eqn:vnm}
\varphi_{a,b}^{op}\equiv \tau_{a,b}\circ \varphi_{a,b}
\end{equation}
where $\tau_{a,b}$ denotes the flip from $A_{a}\otimes A_{b}$
to $A_{b}\otimes A_{a}$.
\item
$\{(A_{a},\varphi_{a,b}):a,b\in\sem\}$
is locally quasi-cocommutative
if there exists $\{R^{(a,b)}:a,b\in\sem\}$
such that $R^{(a,b)}$ is a unitary in
$A_{a}\otimes A_{b}$ 
and
%
%
\begin{equation}
\label{eqn:rabv}
R^{(a,b)}\varphi_{a,b}(x)(R^{(a,b)})^{*}=\varphi_{b,a}^{op}(x)\quad
(x\in A_{ab})
\end{equation}
for each $a,b\in\sem$.
In this case,
we write
$\{(A_{a},\varphi_{a,b},R^{(a,b)}):a,b\in\sem\}$
as a locally  quasi-cocommutative C$^{*}$-WCS. 
\item
A locally  quasi-cocommutative C$^{*}$-WCS 
$\{(A_{a},\varphi_{a,b},R^{(a,b)}):a,b\in\sem\}$
is locally  quasi-triangular
if 
the following holds for each $a,b,c\in\sem$:
%
%
\begin{eqnarray}
\label{eqn:varphiab}
(\varphi_{a,b}\otimes id_{c})(R^{(ab,c)})=&R^{(a,c)}_{13}R^{(b,c)}_{23},\\
\nonumber
\\
\label{eqn:varphiabc}
(id_{a}\otimes \varphi_{b,c})(R^{(a,bc)})=&R^{(a,c)}_{13}R^{(a,b)}_{12}.
\end{eqnarray}
%
\item
A locally  quasi-cocommutative C$^{*}$-WCS 
$\{(A_{a},\varphi_{a,b},R^{(a,b)}):a,b\in\sem\}$
is locally  triangular
if 
$\{(A_{a},\varphi_{a,b},R^{(a,b)}):a,b\in\sem\}$
is locally  quasi-triangular and the following holds:
%
%
\begin{equation}
\label{eqn:tautau}
R^{(a,b)}\tau_{b,a} (R^{(b,a)})=I_{a}\otimes I_{b}\quad(a,b\in\sem)
\end{equation}
where $I_{x}$ denotes the unit of $A_{x}$ for $x=a,b$.
\end{enumerate}
\end{defi}

For a C$^{*}$-WCS 
$\{(A_{a},\varphi_{a,b}):a,b\in\sem\}$,
we see that 
${\cal M}(A_{*}\otimes A_{*})\cong \prod_{a,b\in\sem}A_{a}\otimes A_{b}$
from (\ref{eqn:multiplus}).
Hence we identify an element in ${\cal M}(A_{*}\otimes A_{*})$
with that in $\prod_{a,b\in\sem}A_{a}\otimes A_{b}$.
Then the following holds.
%
%
\begin{lem}(\cite{TS21}, Lemma 2.4)
\label{lem:quasi}
Assume that a monoid $\sem$ is abelian.
\begin{enumerate}
\item
If a C$^{*}$-WCS $\{(A_{a},\varphi_{a,b}):a,b\in\sem\}$
is locally  quasi-cocommutative with respect to $\{R^{(a,b)}:a,b\in\sem\}$ in
(\ref{eqn:rabv}),
then the unitary 
$R\in {\cal M}(A_{*}\otimes A_{*})$ defined by
%
%
\begin{equation}
\label{eqn:sum}
R\equiv (R^{(a,b)})_{a,b\in\sem}
\end{equation}
is a universal $R$-matrix of $(A_{*},\delp)$ in Theorem \ref{Thm:mainthree}.
\item
If a locally  quasi-cocommutative C$^{*}$-WCS 
$\{(A_{a},\varphi_{a,b},R^{(a,b)}):a,b\in\sem\}$
is locally quasi-triangular ({\it resp.} locally triangular),
then 
$(A_{*},\delp,R)$ is quasi-triangular ({\it resp.} triangular)
for $R$ in (\ref{eqn:sum}).
\end{enumerate}
\end{lem}

%
%
\ssft{Componentwise tensor power of C$^{*}$-weakly coassociative system}
\label{subsection:thirdtwo}
In this subsection,
we give a new method to construct
C$^{*}$-weakly coassociative systems (=C$^{*}$-WCSs) 
from a given C$^{*}$-WCS.
Assume that 
$\{(A_{a},\varphi_{a,b}):a,b\in\sem\}$
is a C$^{*}$-WCS.
Fix $n\geq 1$. 
Let $A^{\otimes n}_{a}$ denote the
$n$-times tensor power of $A_{a}$ for $a\in \sem$.
For $a,b\in\sem$,
define $\varphi_{a,b}^{(n)}\in {\rm Hom}(A_{ab}^{\otimes n},
A_{a}^{\otimes n}\otimes A_{b}^{\otimes n})$ by
%
%
\begin{equation}
\label{eqn:vtensor}
\varphi_{a,b}^{(n)}\equiv 
T_{a,b}^{(n)}\circ (\varphi_{a,b})^{\otimes n}
\end{equation}
where $T_{a,b}^{(n)}\in 
{\rm Hom}((A_{a}\otimes A_{b})^{\otimes n}, 
A_{a}^{\otimes n}\otimes A_{b}^{\otimes n})$ 
is defined as 
%
%
\begin{equation}
\label{eqn:taun}
T_{a,b}^{(n)}(x_{1}\otimes y_{1}\otimes x_{2}\otimes y_{2}
\otimes\cdots \otimes
x_{n}\otimes y_{n})
\equiv 
x_{1}\otimes \cdots\otimes x_{n}\otimes 
y_{1}\otimes \cdots\otimes y_{n}
\end{equation}
for $x_{1},\ldots,x_{n}\in A_{a}$ and $y_{1},\ldots,y_{n}\in A_{b}$.
Then we see that 
$(\varphi_{a,b}^{(n)}\otimes id_{c}^{\otimes n})
\circ \varphi_{ab,c}^{(n)}=
(id_{a}^{\otimes n}\otimes \varphi_{b,c}^{(n)})
\circ \varphi_{a,bc}^{(n)}$ for each $a,b,c\in \sem$.
Hence we can verify that
$\{(A_{a}^{\otimes n},\varphi_{a,b}^{(n)}):a,b\in\sem\}$
is a C$^{*}$-WCS.
%
%
\begin{defi}
\label{defi:ctpone}
The C$^{*}$-WCS
$\{(A_{a}^{\otimes n},\varphi_{a,b}^{(n)}):a,b\in\sem\}$
is called the componentwise $n$-times tensor power
of $\{(A_{a},\varphi_{a,b}):a,b\in\sem\}$.
\end{defi}

\noindent
Clearly,
$(\bigoplus_{a} A_{a})^{\otimes n}$ and 
$(\bigoplus_{a}A^{\otimes n}_{a})$ are not isomorphic as a C$^{*}$-algebra
when $n\geq 2$ in general.
Hence the C$^{*}$-bialgebra associated with 
$\{(A_{a}^{\otimes n},\varphi_{a,b}^{(n)}):a,b\in\sem\}$
is not isomorphic to a tensor power of the C$^{*}$-bialgebra 
associated with 
$\{(A_{a},\varphi_{a,b}):a,b\in\sem\}$ in general.
%
%
\begin{lem}
\label{lem:composite}
Assume that 
$\{(A_{a},\varphi_{a,b},R^{(a,b)}):a,b\in\sem\}$
is a locally quasi-cocommutative C$^{*}$-WCS.
Fix $n\geq 1$.
\begin{enumerate}
\item
For $a,b\in\sem$,
$(\varphi_{b,a}^{(n)})^{op}=(\varphi_{b,a}^{op})^{(n)}$.
\item
For $a,b\in\sem$,
define $R^{(a,b:n)}\in A_{a}^{\otimes n}\otimes A_{b}^{\otimes n}$ by
%
%
\begin{equation}
\label{eqn:rabn}
R^{(a,b:n)}\equiv 
T_{a,b}^{(n)}((R^{(a,b)})^{\otimes n}).
\end{equation}
Then $\{(A_{a}^{\otimes n},\varphi_{a,b}^{(n)},R^{(a,b:n)}):
a,b\in\sem\}$ is a locally quasi-cocommutative C$^{*}$-WCS.
\item
In addition to (ii),
if $\{(A_{a},\varphi_{a,b},R^{(a,b)}):a,b\in\sem\}$
is locally quasi-triangular ({\it resp.} locally triangular),
then 
$\{(A_{a}^{\otimes n},\varphi_{a,b}^{(n)},R^{(a,b:n)}):
a,b\in\sem\}$
is also locally quasi-triangular ({\it resp.} locally triangular).
\end{enumerate}
\end{lem}
%
%
\pr
(i)
Let $\tau_{b,a}^{(n)}$ denote 
the flip from
$A_{b}^{\otimes n}\otimes A_{a}^{\otimes n}$ 
to $A_{a}^{\otimes n}\otimes A_{b}^{\otimes n}$.
Then we can verify that
%
%
\begin{equation}
\label{eqn:eone}
T_{a,b}^{(n)}\circ (\tau_{b,a})^{\otimes n}=\tau_{b,a}^{(n)}\circ T_{b,a}^{(n)}.
\end{equation}
On the other hand,
we see that 
$(\varphi_{b,a}^{(n)})^{op}
=
\tau_{b,a}^{(n)}\circ T_{b,a}^{(n)}\circ (\varphi_{b,a})^{\otimes n}$ and
$(\varphi_{b,a}^{op})^{(n)}=
T_{a,b}^{(n)}\circ (\tau_{b,a})^{\otimes n}\circ (\varphi_{b,a})^{\otimes n}$.
From these, the statement holds.

\noindent
(ii)
Let $x=x_{1}\otimes \cdots \otimes x_{n}
\in A_{ab}^{\otimes n}$.
Then \\
$R^{(a,b:n)}\varphi_{a,b}^{(n)}(x)(R^{(a,b:n)})^{*}$
\[
\begin{array}{rl}
=&
T_{a,b}^{(n)}\{(R^{(a,b)})^{\otimes n}(\varphi_{a,b})^{\otimes n}(x)
((R^{(a,b)})^{*})^{\otimes n}\}\\
=&
T_{a,b}^{(n)}\{
R^{(a,b)}\varphi_{a,b}(x_{1})(R^{(a,b)})^{*}
\otimes \cdots \otimes R^{(a,b)}\varphi_{a,b}(x_{n})(R^{(a,b)})^{*}
\}\\
=&
T_{a,b}^{(n)}\{\varphi_{b,a}^{op}(x_{1})
\otimes \cdots\otimes \varphi_{b,a}^{op}(x_{n})\}
\quad(\mbox{from (\ref{eqn:rabv})})
\\
=&(\varphi_{b,a}^{op})^{(n)}(x)\\
=&(\varphi_{b,a}^{(n)})^{op}(x)\quad(\mbox{from (i)}).\\
\end{array}
\]
This implies the statement.

\noindent
(iii)
Assume the local quasi-triangularity for 
$\{(A_{a},\varphi_{a,b},R^{(a,b)}):a,b\in\sem\}$.
For $a,b,c\in\sem$,
define 
$T_{(a,b),c}^{(n)}\in {\rm Hom}((A_{a}\otimes A_{b}\otimes A_{c})^{\otimes n},
(A_{a}\otimes A_{b})^{\otimes n}\otimes A_{c}^{\otimes n})$
and
$T_{a,b,c}^{(n)}\in {\rm Hom}((A_{a}\otimes A_{b}\otimes A_{c})^{\otimes n}, 
A_{a}^{\otimes n}\otimes A_{b}^{\otimes n}\otimes A_{c}^{\otimes n})$
by
%
%
\begin{eqnarray}
&
\label{eqn:tabc}
\begin{array}{l}
T_{(a,b),c}^{(n)}(
x_{1}\otimes y_{1}\otimes z_{1}\otimes \cdots
\otimes 
x_{n}\otimes y_{n}\otimes z_{n})\qquad \\
\qquad \qquad \equiv 
x_{1}\otimes y_{1}\otimes  \cdots \otimes x_{n}\otimes y_{n}
\otimes
z_{1}\otimes \cdots \otimes z_{n},
\end{array}
&
\\
\nonumber
\\
\label{eqn:tabl}
&
\begin{array}{l}
T_{a,b,c}^{(n)}(
x_{1}\otimes y_{1}\otimes z_{1}\otimes \cdots
\otimes 
x_{n}\otimes y_{n}\otimes z_{n})\qquad \\
\qquad \qquad \equiv 
x_{1}\otimes \cdots \otimes x_{n}
\otimes
y_{1}\otimes \cdots \otimes y_{n}
\otimes
z_{1}\otimes \cdots \otimes z_{n}
\end{array}
\end{eqnarray}
for $x_{1},\ldots,x_{n}\in A_{a}$,
$y_{1},\ldots,y_{n}\in A_{b}$
and $z_{1},\ldots,z_{n}\in A_{c}$.
Then \\
$(\varphi_{a,b}^{(n)}\otimes id_{c}^{\otimes n})(R^{(ab,c:n)})$
\[
\begin{array}{rl}
=&
\{(T_{a,b}^{(n)}\otimes id_{c}^{\otimes n})\circ 
((\varphi_{a,b})^{\otimes n}\otimes id_{c}^{\otimes n})\}
(T_{ab,c}^{(n)}((R^{(ab,c)})^{\otimes n}))\\
=&
(T_{a,b}^{(n)}\otimes id_{c}^{\otimes n})
(T_{(a,b),c}^{(n)}(\{(\varphi_{a,b}\otimes id_{c})(R^{(ab,c)})\}^{\otimes n}))\\
=&
(T_{a,b}^{(n)}\otimes id_{c}^{\otimes n})
(T_{(a,b),c}^{(n)}(\{R_{13}^{(a,c)}R_{23}^{(b,c)}
\}^{\otimes n}))\quad(\mbox{from (\ref{eqn:varphiab})})\\
=&
T_{a,b,c}^{(n)}(\{R_{13}^{(a,c)}R_{23}^{(b,c)}\}^{\otimes n})\\
=&
T_{a,b,c}^{(n)}(\{R_{13}^{(a,c)}\}^{\otimes n})
T_{a,b,c}^{(n)}(\{R_{23}^{(b,c)}\}^{\otimes n})\\
=&
R^{(a,c:n)}_{13}R^{(b,c:n)}_{23}.
\end{array}
\]
By the same reasoning,
we obtain 
$(id_{a}^{\otimes n}\otimes \varphi_{b,c}^{(n)})(R^{(a,bc:n)})
=R^{(a,c:n)}_{13}R^{(a,b:n)}_{12}$.
Hence the statement about the local quasi-triangularity holds.

Assume the local triangularity for 
$\{(A_{a},\varphi_{a,b},R^{(a,b)}):a,b\in\sem\}$.
It is sufficient to show (\ref{eqn:tautau}) for
$\{(A_{a}^{\otimes n},\varphi_{a,b}^{(n)},R^{(a,b:n)}):a,b\in\sem\}$
here.
For $a,b\in\sem$,
let $\tau_{b,a}^{(n)}$ be as in the proof of (i).
Then\\
$R^{(a,b:n)}\tau_{b,a}^{(n)}(R^{(b,a:n)})$
\[
\begin{array}{rl}
=&
T_{a,b}^{(n)}((R^{(a,b)})^{\otimes n})
\,\tau_{b,a}^{(n)}\{T_{b,a}^{(n)}((R^{(b,a)})^{\otimes n})\}\\
=&
T_{a,b}^{(n)}((R^{(a,b)})^{\otimes n})\,
(T_{a,b}^{(n)}\circ (\tau_{b,a})^{\otimes n})\{(R^{(b,a)})^{\otimes n}\}\quad
(\mbox{from (\ref{eqn:eone})})\\
=&
T_{a,b}^{(n)}\{
(R^{(a,b)})^{\otimes n}(\tau_{b,a})^{\otimes n}((R^{(b,a)})^{\otimes n})\}
\\
=&
T_{a,b}^{(n)}(\{R^{(a,b)}\tau_{b,a}(R^{(b,a)})\}^{\otimes n})\\
=&
T_{a,b}^{(n)}(\{I_{a}\otimes I_{b}\}^{\otimes n})\quad(\mbox{from (\ref{eqn:tautau})})\\
=& 
I_a^{\otimes n}\otimes I_b^{\otimes n}.
\end{array}
\]
Hence the statement about the local triangularity holds.
\qedh

%
%
\ssft{Componentwise infinite tensor power of C$^{*}$-weakly coassociative system}
\label{subsection:thirdthree}
In this subsection,
we define the componentwise infinite tensor power of C$^{*}$-weakly coassociative system
(=C$^{*}$-WCS) from a given C$^{*}$-WCS
as the inductive limit of 
componentwise tensor powers of C$^{*}$-WCS.
Let  
$\{(A_{a},\varphi_{a,b}):a,b\in\sem\}$ be a C$^{*}$-WCS and 
let $\{(A_{a}^{\otimes n},\varphi_{a,b}^{(n)}):a,b\in\sem\}$
be the componentwise $n$-times tensor power
of $\{(A_{a},\varphi_{a,b}):a,b\in\sem\}$ in Definition \ref{defi:ctpone}
for $n\geq 1$.
With respect to the embedding 
%
%
\begin{equation}
\label{eqn:psib}
\psi^{(n)}_{a}:A_{a}^{\otimes n}\ni x\mapsto
 x\otimes I_{a}\in A_{a}^{\otimes n}\otimes A_{a}
=A_{a}^{\otimes (n+1)},
\end{equation}
we regard $A_{a}^{\otimes n}$ as a C$^{*}$-subalgebra of 
$A_{a}^{\otimes (n+1)}$ for each $a\in\sem$.
Let $A_{a}^{\otimes \infty}$ denote the inductive limit of the inductive system
$\{(A_{a}^{\otimes n},\psi_{a}^{(n)}):n\geq 1\}$ of C$^{*}$-algebras:
%
%
\begin{equation}
\label{eqn:abb}
A_a^{\otimes \infty}\equiv  \varinjlim 
(A_{a}^{\otimes n},\psi_{a}^{(n)}).
\end{equation}
The C$^{*}$-algebra $A_a^{\otimes \infty}$ is called the 
{\it infinite tensor product} of $A_{a}$ (\cite{Blackadar2006}, $\S$ II.9.8).
The map $\psi_{a}^{(n)}$ in (\ref{eqn:psib}) satisfies 
%
%
\begin{equation}
\label{eqn:psitwob}
(\psi_{a}^{(n)}\otimes \psi_{b}^{(n)})\circ 
\varphi^{(n)}_{a,b}
=
\varphi^{(n+1)}_{a,b}\circ \psi_{ab}^{(n)}\quad
(a,b\in\sem,\,n\geq 1).
\end{equation}
From $\{\varphi_{a,b}^{(n)}:n\geq 1\}$,
we can define the $*$-homomorphism $\varphi_{a,b}^{(\infty)}$ 
from $A_{ab}^{\otimes \infty}$
to $A_{a}^{\otimes \infty}\otimes A_{b}^{\otimes \infty}$ 
such that 
%
%
\begin{equation}
\label{eqn:isomorphismb}
(\varphi_{a,b}^{(\infty)})|_{A_{ab}^{\otimes n}}
=\varphi_{a,b}^{(n)}
\end{equation}
for each $n$
where we identify 
$A_{a}^{\otimes \infty}\otimes A_{b}^{\otimes \infty}$ 
with the inductive limit of the system
$\{(A_{a}^{\otimes n}\otimes A_{b}^{\otimes n},
\psi^{(n)}_{a}\otimes \psi^{(n)}_{b}):n\geq 1\}$.
Then the following holds:
%
%
\begin{equation}
\label{eqn:commutativeb}
(\varphi_{a,b}^{(\infty)}\otimes id_{c})\circ \varphi_{ab,c}^{(\infty)}
=
(id_{a}\otimes \varphi_{b,c}^{(\infty)})\circ \varphi_{a,bc}^{(\infty)}
\quad(a,b,c\in\sem)
\end{equation}
where $id_{x}$ denotes the identity map 
on $A_{x}^{\otimes \infty}$ for $x=a,c$.
From this,
we see that 
$\{(A_a^{\otimes \infty},\varphi^{(\infty)}_{a,b}):a,b\in\sem\}$
is a C$^{*}$-WCS.
%
%
\begin{defi}
\label{defi:ctptwo}
The C$^{*}$-WCS
$\{(A_{a}^{\otimes \infty},\varphi_{a,b}^{(\infty)}):a,b\in\sem\}$
is called the componentwise infinite tensor power
of $\{(A_{a},\varphi_{a,b}):a,b\in\sem\}$.
\end{defi}

By Theorem \ref{Thm:mainthree}, the following direct sum
%
%
\begin{equation}
\label{eqn:sumsumb}
(A^{\otimes \infty})_*\equiv \bigoplus_{a\in \sem}A_a^{\otimes \infty}
\end{equation}
is a C$^{*}$-bialgebra.
From here,
we write $A^{\otimes \infty}_*$ as
$(A^{\otimes \infty})_*$ for simplicity of description.

Let $\psi_{*}^{(n)}\equiv \bigoplus_{a\in \sem}\psi_{a}^{(n)}$
and $A^{\otimes n}_*\equiv \bigoplus_{a\in\sem}A^{\otimes n}_a$.
Then $\{(A^{\otimes n}_*,\psi_{*}^{(n)}):n\geq 1\}$ is 
an inductive system of C$^{*}$-bialgebras.
We see that
$A^{\otimes \infty}_*$ in (\ref{eqn:sumsumb})
is the inductive limit of 
the inductive system
$\{(A_{*}^{\otimes n},\psi_{*}^{(n)}):n\geq 1\}$ of C$^{*}$-algebras.
Furthermore,
we can verify that
the C$^{*}$-bialgebra associated with the C$^{*}$-WCS
$\{(A_a^{\otimes \infty},\varphi^{(\infty)}_{a,b}):a,b\in\sem\}$
coincides
with the inductive limit $\varinjlim (A^{\otimes n}_*,\psi_{*}^{(n)})$
of C$^{*}$-bialgebras
$\{(A^{\otimes n}_*,\psi_{*}^{(n)}):n\geq 1\}$.
Hence we obtain the following equation of C$^{*}$-bialgebras:
%
%
\begin{equation}
\label{eqn:last}
A^{\otimes \infty}_{*}=
\varinjlim (A^{\otimes n}_*,\psi_{*}^{(n)}).
\end{equation}
%
%
%
\sftt{Proof of Theorem \ref{Thm:mainone}}
\label{section:fourth}
In this section,
we prove Theorem \ref{Thm:mainone}.
Procedures are as follows:
In $\S$ \ref{subsection:fourthone},
we recall the example of locally triangular 
C$^{*}$-weakly coassociative system (=C$^{*}$-WCS) in \cite{TS21}.
From this, 
an inductive system of triangular (especially, quasi-cocommutative) 
C$^{*}$-bialgebras is constructed
by using the method in $\S$ \ref{subsection:thirdthree}.
In $\S$ \ref{subsection:fourthtwo},
it is proved that its inductive limit is not quasi-cocommutative.
%
%
%
\ssft{An example of locally triangular C$^{*}$-weakly coassociative system}
\label{subsection:fourthone}
Let ${\bf N}\equiv \{1,2,3,\ldots\}$.
We regard ${\bf N}$ as an abelian monoid with respect to the multiplication.
For $n\in {\bf N}$,
let $M_{n}$ denote the (finite-dimensional)
C$^{*}$-algebra of all $n\times n$-complex matrices
where we define $M_{1}={\bf C}$.
We recall a locally triangular C$^{*}$-WCS 
$\{(M_{n}, \varphi_{n,m},R^{(n,m)}):n,m\in {\bf N}\}$
in \cite{TS21} as follows.

Let $\{E^{(n)}_{i,j}\}_{i,j=1}^{n	}$ denote
the set of standard matrix units of $M_{n}$.
For $n,m\in {\bf N}$,
define the $*$-isomorphism $\varphi_{n,m}$ 
from $M_{nm}$ onto $M_{n}\otimes M_{m}$ by
%
%
\begin{equation}
\label{eqn:eij}
\varphi_{n,m}(E_{m(i-1)+j,m(i^{'}-1)+j^{'}}^{(nm)})=
E_{i,i^{'}}^{(n)}\otimes  E_{j,j^{'}}^{(m)}
\end{equation}
for $i,i^{'}\in\{1,\ldots,n\}$ and $j,j^{'}\in\{1,\ldots,m\}$.
For $n\in {\bf N}$, let $\{e^{(n)}_{i}\}_{i=1}^{n}$ denote 
the standard basis of the finite dimensional Hilbert space ${\bf C}^{n}$.
Define the unitary transformation $R^{(n,m)}$ on 
${\bf C}^{n}\otimes {\bf C}^{m}$ by
%
%
\begin{equation}
\label{eqn:urm}
R^{(n,m)}(e_{i}^{(n)}\otimes e_{j}^{(m)})\equiv 
e_{\ul{i}}^{(n)}\otimes e_{\ul{j}}^{(m)}
\end{equation}
for $(i,j)\in \{1,\ldots,n\}\times\{1,\ldots,m\}$
where the pair $(\ul{i},\ul{j})\in \{1,\ldots,n\}\times
\{1,\ldots,m\}$ is uniquely defined as the following integer equation:
%
%
\begin{equation}
\label{eqn:mij}
m(i-1)+j=n(\ul{j}-1)+\ul{i}.
\end{equation}
By the natural identification 
${\rm End}_{{\bf C}}({\bf C}^{n}\otimes {\bf C}^{m})
\cong M_{n}\otimes M_{m}$,
$R^{(n,m)}$ is regarded as 
a unitary element in  $M_{n}\otimes M_{m}$ for each $n,m\in {\bf N}$.
Then $\{(M_{n},\varphi_{n,m},R^{(n,m)}):n,m\in {\bf N}\}$
is a locally triangular WCS (\cite{TS21}, $\S$ 3).

For $i\geq 1$, let $M_{n}^{\otimes i}$ denote the $i$-times tensor power of $M_{n}$.
From Lemma \ref{lem:composite}(iii),
we obtain the locally triangular C$^{*}$-WCS 
 $\{(M_{n}^{\otimes i},\varphi^{(i)}_{n,m},R^{(n,m:i)}):n,m\in {\bf N}\}$
associated with 
$\{(M_{n},\varphi_{n,m},R^{(n,m)}):n,m\in {\bf N}\}$.
By Lemma \ref{lem:quasi}(ii), 
%
%
\begin{equation}
\label{eqn:istar}
M^{\otimes i}(*)\equiv \bigoplus_{n\in {\bf N}}M_n^{\otimes i}
\end{equation}
is the triangular C$^{*}$-bialgebra associated with 
$\{(M_{n}^{\otimes i},\varphi_{n,m}^{(i)},R^{(n,m:i)}):n,m\in {\bf N}\}$.
Especially, 
$M^{\otimes i}(*)$ is quasi-cocommutative
for each $1\leq i<\infty$.

Define the C$^{*}$-algebra $M^{\otimes \infty}(*)$ by
%
%
\begin{equation}
\label{eqn:sumsum}
M^{\otimes \infty}(*)\equiv \bigoplus_{n\in {\bf N}}M_n^{\otimes \infty}
\end{equation}
where $M_n^{\otimes \infty}$ denotes the infinite tensor product of 
the C$^{*}$-algebra $M_n$.
By definition, $M_{n}^{\otimes \infty}$ is 
a uniformly hyperfinite algebra of Glimm's type $\{n^{i}\}_{i\geq 1}$ 
\cite{Glimm,TS18}.
From $\S$ \ref{subsection:thirdthree},
$M^{\otimes \infty}(*)$ is the inductive limit of 
the inductive system $\{(M^{\otimes i}(*),\psi_{*}^{(i)}):i\geq 1\}$
of C$^{*}$-bialgebras: 
%
%
\begin{equation}
\label{eqn:id}
M^{\otimes \infty}(*)=
\varinjlim (M^{\otimes i}(*),\psi_{*}^{(i)}).
\end{equation}

We illustrate relations among algebras as follows:

\def\inclusion{
\put(-30,500){$M(*)$}
\put(-30,400){$M^{\otimes 2}(*)$}
\put(-30,300){$M^{\otimes 3}(*)$}
\put(-10,200){$\vdots$}
\put(-30,100){$M^{\otimes \infty}(*)$}
\put(10,480){\rotatebox{-90}{$\hookrightarrow$}}
\put(10,380){\rotatebox{-90}{$\hookrightarrow$}}
\put(10,280){\rotatebox{-90}{$\hookrightarrow$}}
\put(120,500){$=$}
\put(120,400){$=$}
\put(120,300){$=$}
\put(120,100){$=$}
\put(170,500){$M_{1}$}
\put(170,400){$M_{1}^{\otimes 2}$}
\put(170,300){$M_{1}^{\otimes 3}$}
\put(180,200){$\vdots$}
\put(160,100){$M_{1}^{\otimes \infty}$}
\put(180,480){\rotatebox{-90}{$=$}}
\put(180,380){\rotatebox{-90}{$=$}}
\put(180,280){\rotatebox{-90}{$=$}}
\put(280,500){$\oplus$}
\put(280,400){$\oplus$}
\put(280,300){$\oplus$}
\put(280,100){$\oplus$}
\put(330,500){$M_{2}$}
\put(330,400){$M_{2}^{\otimes 2}$}
\put(330,300){$M_{2}^{\otimes 3}$}
\put(340,200){$\vdots$}
\put(320,100){$M_{2}^{\otimes \infty}$}
\put(340,480){\rotatebox{-90}{$\hookrightarrow$}}
\put(340,380){\rotatebox{-90}{$\hookrightarrow$}}
\put(340,280){\rotatebox{-90}{$\hookrightarrow$}}
\put(440,500){$\oplus$}
\put(440,400){$\oplus$}
\put(440,300){$\oplus$}
\put(440,100){$\oplus$}
\put(500,500){$M_{3}$}
\put(500,400){$M_{3}^{\otimes 2}$}
\put(500,300){$M_{3}^{\otimes 3}$}
\put(520,200){$\vdots$}
\put(490,100){$M_{3}^{\otimes \infty}$}
\put(510,480){\rotatebox{-90}{$\hookrightarrow$}}
\put(510,380){\rotatebox{-90}{$\hookrightarrow$}}
\put(510,280){\rotatebox{-90}{$\hookrightarrow$}}
\put(600,500){$\oplus \cdots$}
\put(600,400){$\oplus \cdots$}
\put(600,300){$\oplus\cdots$}
\put(600,100){$\oplus\cdots$}
%
\thicklines
\put(-180,450){\vector(0,-1){320}}
\put(-210,520){{\small {\sf inductive}}}
\put(-210,480){{\small {\sf limit}}}
}
\thicklines
\setlength{\unitlength}{.1mm}
\begin{picture}(1000,550)(-300,50)
\put(0,0){\inclusion}
\end{picture}

%
%
\ssft{Proof of Theorem \ref{Thm:mainone}}
\label{subsection:fourthtwo}
In this subsection, we prove Theorem \ref{Thm:mainone}.
%
%
\begin{lem}
\label{lem:noq}
Let $\{(M_n^{\otimes \infty},\varphi^{(\infty)}_{n,m}):n,m\in {\bf N}\}$	
 be 
the componentwise infinite tensor power of 
$\{(M_n,\varphi_{n,m}):n,m\in {\bf N}\}$	
in $\S$ \ref{subsection:fourthone}.
\begin{enumerate}
\item
For two representations $\pi_{1}$ and $\pi_{2}$ of 
the C$^{*}$-algebra $M_2^{\otimes \infty}$,
define the representation of $M_{4}^{\otimes \infty}$ by
%
%
\begin{equation}
\label{eqn:sstar}
\pi_{i}\star \pi_{j}\equiv (\pi_{i}\otimes \pi_{j})\circ \varphi_{2,2}^{(\infty)}
\quad(i,j=1,2).
\end{equation}
Then there exist two unital representations 
$\pi_{1}$ and $\pi_{2}$ of $M_2^{\otimes \infty}$ 
such that 
$\pi_{1}\star\pi_{2}$ and 
$\pi_{2}\star\pi_{1}$ are not unitarily equivalent.
\item
Let $\delp$ denote the comultiplication of $M^{\otimes \infty}(*)$ 
in (\ref{eqn:sumsum})
with respect to the C$^{*}$-weakly coassociative system
$\{(M_n^{\otimes \infty},\varphi^{(\infty)}_{n,m}):n,m\in {\bf N}\}$.
Then the C$^{*}$-bialgebra $(M^{\otimes \infty}(*),\delp)$
is not quasi-cocommutative. 
\end{enumerate}
\end{lem}
%
%
\pr
(i)
For $i=1,2$,
define the (pure) state $\omega_{i}^{(0)}$ of $M_{2}$ by
%
%
\begin{equation}
\label{eqn:etwo}
\omega_{i}^{(0)}(x)\equiv x_{ii}
\quad(x=(x_{ij})_{i,j=1}^{2}\in M_{2})
\end{equation}
where $x_{ij}$'s denote standard matrix units of 
the $2\times 2$ matrix $x$.
Let $\omega_{i}$ denote the product state $(\omega_{i}^{(0)})^{\otimes \infty}$
of $M_{2}^{\otimes \infty}$ for $i=1,2$.
Let $\pi_{i}$ denote the Gel'fand-Na\v{\i}mark-Segal
representation of $M_{2}^{\otimes \infty}$ by $\omega_{i}$ and
let $P[i]$ denote its unitary equivalence class.
Then $P_{2}[1]\star P_{2}[2]\ne P_{2}[2]\star P_{2}[1]$  
from  (2.6) and (3.2) in \cite{TS18}
where we remark that $\star$ is well-defined on unitary equivalence classes
of representations.
Hence the statement is proved.

\noindent
(ii)
Let $p_{n}$ denote the projection from $M^{\otimes \infty}(*)$
to $M^{\otimes \infty}_{n}$ for $n\in {\bf N}$.
From this, any representation of $M^{\otimes \infty}_{n}$ lifts on
$M^{\otimes \infty}(*)$.
Let $\pi_{1}$ and $\pi_{2}$ be unital representations of $M_{2}^{\otimes \infty}$
in (i).
Then $(\pi_{i}\star \pi_{j})\circ p_{4}$ is a nondegenerate representation 
of $M^{\otimes \infty}(*)$
such that $(\pi_{i}\star \pi_{j})\circ p_{4}
=(\pi_{i}\circ p_{2}\otimes \pi_{j}\circ p_{2})\circ \delp$.
From this and (i),
$(\pi_{1}\circ p_{2}\otimes \pi_{2}\circ p_{2})\circ \delp$
and 
$(\pi_{2}\circ p_{2}\otimes \pi_{1}\circ p_{2})\circ \delp$
are not unitarily equivalent.
From this and Lemma \ref{lem:maintwob}, 
$(M^{\otimes \infty}(*),\delp)$ is not quasi-cocommutative.
\qedh
%
%
\\
{\it Proof of Theorem \ref{Thm:mainone}.}
By Remark \ref{rem:cat},
the category of (quasi-cocommutative) C$^{*}$-bialgebras
makes sense.
From this and Lemma \ref{lem:noq}(ii) and (\ref{eqn:id}),  
$\{M^{\otimes i}(*)\}_{i\geq 1}$
is an example of inductive system of quasi-cocommutative
C$^{*}$-bialgebras
such that its inductive limit is not quasi-cocommutative.
This example implies the statement.
\qedh

\noindent
From Remark \ref{rem:cat},	
Lemma \ref{lem:noq}(ii) and (\ref{eqn:id}),  
the following is automatically proved.
%
%
\begin{cor}
\label{cor:option}
\begin{enumerate}
\item
The category of quasi-triangular C$^{*}$-bialgebras
is not closed with respect to the inductive limit.
\item
The category of triangular C$^{*}$-bialgebras
is not closed with respect to the inductive limit.
\end{enumerate}
\end{cor}



%
%

%
\end{document}